\documentclass[10pt]{amsart}
\usepackage{a4,xy,amssymb,times}
\xyoption{all}
\def\r{\rightarrow}

\def\Id{\operatorname{id}}

\let\cal\mathcal
\def\Ascr{{\cal A}}

\def\Cscr{{\cal C}}
\def\Dscr{{\cal D}}

\def\Sscr{{\cal S}}
\def\Tscr{{\cal T}}

\newcommand{\se}[1]{\begin{equation*}\begin{split}#1\end{split}\end{equation*}}

\newcommand{\C}{\mathbb{C}}

\newcommand{\atimes}{\underset{\text{{\it\tiny A--A}}}{\otimes}}
\newcommand{\stimes}{\underset{\text{{\it\tiny S--S}}}{\otimes}}
\newcommand{\cD}{\mathcal{D}}
\newcommand{\cI}{\mathcal{I}}
\newcommand{\cJ}{\mathcal{J}}
\newcommand{\cC}{\mathcal{C}}
\newcommand{\cK}{\mathcal{K}}
\newcommand{\cP}{\mathcal{P}}

\newcommand{\vtx}[1]{*+[o][F-]{\scriptscriptstyle #1}}

\newcommand{\Tr}{\textrm{Tr}}

\newcommand{\cR}{\mathcal{R}}
\newcommand{\cycle}{\circlearrowright}
\newcommand{\cycl}{\cycle}

\newcommand{\<}{\langle}
\renewcommand{\>}{\rangle}

\newcommand{\Mod}{\ensuremath{\mathsf{Mod}}}
\newcommand{\Bimod}{\ensuremath{\mathsf{Mod}}}

\newcommand{\Ext}{\mathsf{Ext}}
\newcommand{\Dim}{\mathsf{Dim \,}}

\newtheorem{lemma}{Lemma}[section]

\newtheorem{theorem}[lemma]{Theorem}
\newtheorem{corollary}[lemma]{Corollary}
\newtheorem{property}[lemma]{Property}

\newtheorem{lemmas}{Lemma}[subsection]
\newtheorem{propositions}[lemmas]{Proposition}
\newtheorem{theorems}[lemmas]{Theorem}

\theoremstyle{definition}

\newtheorem{definition}[lemma]{Definition}

\newtheorem{examples}[lemmas]{Example}
\newtheorem{definitions}[lemmas]{Definition}

{

}

\theoremstyle{remark}

\newtheorem{remark}[lemma]{Remark}
\newtheorem{remarks}[lemmas]{Remark}

\newcommand{\Sup}{\ensuremath{\mathsf{Sup}}}

\newcommand{\Mat}{\texttt{Mat}}

\newcommand{\Rep}{\texttt{Rep}}
\newcommand{\RRep}{\texttt{RRep}}
\newcommand{\Coh}{\texttt{Coh}}
\newcommand{\Hom}{\textrm{Hom}}
\newcommand{\RHom}{\textrm{RHom}}

\newcommand{\id}{\mathbf{1}}

\newcommand{\lt}{\ensuremath{\mathsf{lt}}}

\title{Graded Calabi Yau Algebras of dimension $3$}

\author{Raf Bocklandt\\With an appendix by Michel Van den Bergh}
\address{Raf Bocklandt\\University of Antwerp\\Middelheimlaan 1\\B-2020 Antwerpen (Belgium)}
\email{rafael.bocklandt@ua.ac.be}
\thanks{The author is a Postdoctoral Fellow of the Fund for Scientific Research - Flanders (Belgium)}

\xyoption{all}
\begin{document}

\begin{abstract}
In this paper we prove that Graded Calabi Yau Algebras of dimension $3$ are isomorphic to 
path algebras of quivers with relations derived from a superpotential. We show that for a given
quiver $Q$ and a degree $d$, the set of good superpotentials of degree $d$, i.e. those that give rise to Calabi Yau algebras is either empty or almost everything (in the measure theoretic sense). We also give some constraints on the structure of quivers that allow good superpotentials, and for the simplest quivers we give a complete list of the degrees for which good superpotentials exist.
\end{abstract}

\maketitle

\section{Introduction and Motivation}

If one studies boundary conditions of the $B$-model in super string theory over an $n$-dimensional Calabi Yau manifold $X$, one obtains naturally the derived category of coherent sheaves $\cD^b\Coh X$ \cite{orlov}. This category is called a Calabi Yau category of dimension three, i.e. the third shift in the derived category is a Serre Functor: 
\[
\forall A,B \in \cD^b\Coh X:~\Hom_{\cD^b\Coh X}(A,B) \cong \Hom_{\cD^b\Coh X}(B,A[3])^*,
\]
where the isomorphisms are natural in $A$ and $B$.
In general this category is too big to study its structure directly and therefore
it is interesting to look at full triangulated subcategories of $\cD^b\Coh X$ that can be modeled using derived categories of module categories of noncommutative algebras.
In string theoretical papers this is often done using path algebras of quivers with relations coming form a
superpotential: if $Q$ is a quiver and $\C Q$ the corresponding path algebra, then a superpotential is an element of the vector space $\C Q / [\C Q,\C Q]$. On this space we can define for every arrow $a$ a 'derivation' $\partial_a$ that cuts out $a$ (for a precise definition see section \ref{quiverdef}).
Given a superpotential $W$ one can construct the \emph{vacualgebra} \cite{BLS}
\[
A_W := \C Q / (\partial_a W: a \in Q_1).
\]
In the exemplary cases worked out by physicists, the derived category of finite dimensional modules of the vacualgebra is indeed a Calabi Yau category, and hence these algebras are called Calabi Yau Algebras. 

In this note we will show that in the case of graded algebras, every graded path algebra with relations that is Calabi Yau of dimension $3$ must be isomorphic to a vacualgebra of some superpotential. The converse is not true but we will show that being a Calabi Yau algebra of dimension $3$ corresponds to the exactness of a certain bimodule complex. Therefore, for a given quiver $Q$ and a given degree $d$ the subset of superpotentials of degree $d$ that give rise to Calabi Yau vacualgebras is either empty or almost everything. Furthermore we will use Groebner basis techniques to explicitly determine the list of degree of good superpotentials of simple quivers.

The results in this paper build further on ideas introduced by M. Van den Bergh in \cite{oberwolf}.
Similar results on Calabi Yau algebras in different settings have been obtained by R. Rouquier and V. Ginzberg \cite{RG}. 
\section{Preliminaries} 

\subsection{Path Algebras with relations}\label{quiverdef}
As usual a \emph{quiver} $Q$ is an oriented graph. We denote the set of vertices by $Q_0$, the set of arrows by $Q_1$ and the maps $h,t$ assign to each arrow its head and tail.
A \emph{nontrivial path} $p$ is a sequence of arrows $a_1\cdots a_k$ such that $t(a_i)=h(a_{i+1})$, whereas a \emph{trivial path} is just a vertex. We will denote the length of a path by $|p|:= k$ and the head and tail by $h(p)=h(a_1),~ t(p)=t(a_k)$. A path is called a cycle if $h(p)=t(p)$. A quiver is called
\emph{connected} if it is not the disjoint union of two subquivers and it is \emph{strongly connected} if there is a cycle through each pair of vertices.

The path algebra $\C Q$ is the complex vector space with as basis the paths in $Q$ and the multiplication of two paths $p$, $q$ is their concatenation $pq$ if $t(p)=h(q)$ or else $0$.
We can put a gradation on $\C Q$ using the length of the paths. The space spanned by all paths of nonzero length is a graded ideal of $\C Q$ and we will denote it by $\cJ$.

The vector space $\C Q/[\C Q,\C Q]$ has as basis the set of cycles up to cyclic permutation of the arrows. We can embed this space into $\C Q$ by mapping a cycle onto the sum of all its possible cyclic permutations:
\[
\cycle : \C Q/[\C Q,\C Q] \to \C Q: a_1\cdots a_n \mapsto \sum_i a_i\cdots a_na_1\cdots a_{i-1}.
\] 
Another convention we will use is the inverse of arrows: if $p:= a_1\cdots a_n$ is a path and $b$ an arrow, then
$pb^{-1}=a_1\cdots a_{n-1}$ if $b=a_n$ and zero otherwise. Similarly one can define $b^{-1}p$. These new defined maps can be combined to obtain a 'derivation' 
\[
\partial_a : \C Q/[\C Q,\C Q] \to \C Q : p \mapsto \cycle(p)a^{-1} =a^{-1}\cycle(p).
\]

From now on $A$ will denote the quotient algebra $\C Q/\cI$ by
a finitely generated graded ideal $\cI \subset \cJ^2$.  
The set $\cR \subset \cI$ will be a minimal set of homogeneous generators each sitting inside some $i\C Qj,~i,j \in Q_0$. 

We denote the semi-simple (left) $A$-module $A/A_{\ge 1}\cong \C Q/\cJ$ by $S$. $S$ is the direct sum of $\#Q_0$ simple one-dimensional $A$-modules $S_i$, each corresponding to a vertex $i \in \C Q$. To each vertex we can also assign a projective module $P_i$ which is the left ideal $Ai$ and $S_i = P_i/{(P_i)}_{\ge 1}$. Although it is a little sloppy we will also use $S$ to denote the subring $A_0\cong \C Q_0$, generated by the vertices.
 
\subsection{Calabi Yau Categories}

Let $\cC$ be an abelian $\C$-linear category and $\cD^b \cC$ its bounded derived category. 
Using the shift we can define a graded functor $(s, \eta^s)$ in the sense of \ref{appendix} where $s$ is the shift functor and the $\eta^s$ gives natural isomorphisms
\[
\eta^s_A : s(A[1]) \to (sA)[1]: x \mapsto -x. 
\]
As explained in the appendix, these maps are uniquely determined by the
demand of compatibility with the triangulated structure of $\cD^b \cC$.

\begin{definition}
The category $\cD^b \cC$ is called \emph{Calabi Yau of dimension $n$} if 
there are natural isomorphisms 
\[
\nu_{A,B}: \Hom_{\cD^b \cC}(A,B) \to \Hom_{\cD^b \cC}(B,s^n A)^*, \text{ ($^*$ is the complex dual)}
\]
or, in other words, the $n^{th}$ shift is a Serre Functor.
\end{definition}

Starting with a graded path algebra with relations $A$, we can construct the category of finite dimensional left $A$-modules: $\Rep A$. This is an abelian category so we can construct its bounded derived category $\cD^b\Rep A$. We will call $A$ a graded Calabi Yau Algebra of dimension $n$ if $\cD^b\Rep A$ is a Calabi Yau category of dimension $n$.

Although the definition is asymmetric in the sense that one only uses left modules, it is easy to see
that if $A$ is Calabi Yau, the derived category of finite dimensional \emph{right} modules $\cD^b \RRep A$ 
is also a Calabi Yau Category. This can be proved using the complex dual as an anti-equivalence between
$\cD^b\Rep A$ and $\cD^b \RRep A$: let $M,N$ be complexes of right modules and define 
\[
\nu^{\RRep A}_{M,N}: \Hom_{\cD^b \RRep A}(M,N) \to \Hom_{\cD^b \RRep A}(N,s^n M)^*
\]
by the equality
\[
\nu^{\RRep A}_{M,N}(f)(g) =\nu^{\Rep A}_{N^*,M^*}(s^n(f^T))(s^n(g^T)).
\]

The Calabi Yau property of the derived category can be tracked back to the original category 
to give us properties that we will often use
\begin{property}
If $A$ is Calabi Yau of dimension $n$ then 
\begin{itemize}
\item[C1]
The global dimension of $A$ is also $n$.
\item[C2] 
If $X, Y \in \Rep A$ then
\[
\Ext^k_A(X,Y) \cong \Ext^{n-k}_A(Y,X)^*.
\]
\item[C3]
The identifications above gives us a pairings $\<,\>^k_{XY}: \Ext^k_A(X,Y) \times \Ext^{n-k}_A(Y,X) \to \C$ which satisfy
\[
\<f,g\>^k_{XY} = \<1_X,g*f\>^0_{XX} = (-1)^{k(n-k)}\<1_Y,f*g\>^0_{YY}, 
\]
where $*$ denotes the standard composition of extensions.
\end{itemize}
\end{property}
\begin{proof}
$(1):$ if $i>n$ then $\Ext^{i}_A(M,N)=\Ext^{n-i}(M,N)=0$ so $\rm{gldim}A\le n$ and 
$\Ext^{n}_A(A/A_+,A/A_+)=\Hom_A(A/A_+,A/A_+)=A/A_+\ne 0$ so $\rm{gldim}A\ge n$.
For $(2-3)$ see the appendix.
\end{proof}

\section{Graded Calabi Yau algebra's of dimension $n \le 3$}
In this section we will give descriptions of the types of quivers and relations that appear in graded Calabi Yau algebras of dimension $3$.

From now on we will also assume that the quiver $Q$ is connected. This is not a severe restriction because $A$ is the direct sum of subalgebras defined over its connected components. Many properties like the Calabi Yau property transfer from the algebra to its direct summands: $A_1\oplus A_2$ is Calabi Yau of dimension $n$ if both $A_1$ and $A_2$ are Calabi Yau of dimension $n$. This follows from the fact
that the representation category (and hence the derived category) of $A$ 
decomposes as the direct sum of $\Rep A_1$ and $\Rep A_2$.

\begin{theorem}\label{nodige}
If $A$ is Calabi Yau of dimension $3$ then 
\begin{enumerate}
\item
there is a homogeneous superpotential $W \in \C Q/[\C Q,\C Q]$ such that 
\[
A \cong \C Q/(\partial_a W: a\in Q_1), 
\]
\item
every arrow in $Q$ is contained in a cycle of $\cycle W$.
\item
every vertex in $Q$ is the source of two arrows and the target of two arrows.
\end{enumerate}
\end{theorem}
\begin{proof}
As the global dimension of $A$ must be $3$, there is a projective graded resolution
\[
\xymatrix{
\bigoplus_{j \in Q_0}P_j^{m_{ij}}\ar@{^(->}[r]^{(f_r)}&\bigoplus_{t(r)=i}P_{h(r)}\ar[r]^{(rb^{-1})}&\bigoplus_{t(b)=i}P_{h(b)}\ar[r]^{(\cdot b)}&P_i\ar@{->>}[r]&S_i.
}
\]
In the diagram above the $r's$ are the relations in $\cR$ and the $b's$ are arrows, the $f_r$ are maps that are not further specified. 
Using the Calabi Yau property and comparing dimensions we can conclude that
\begin{enumerate}
\item $m_{ij}=\Dim \Ext^3(S_i,S_j) = \Dim \Hom(S_j,S_i) = \delta_{ij}$,
\item $\#\{r\in \cR: h(r)=j,t(r)=i \}=\Dim \Ext^2(S_i,S_j)=\Dim \Ext^1(S_i,S_j)= \#\{a\in Q_1:i \stackrel{a}{\leftarrow} j\}$.
\end{enumerate}
Because of $(1)$ we can identify each $f_r$ with an element in $i A h(r)$. 
Consider the finite dimensional quotient algebra 
\[
M = A/(f_r: r \in \cR, A_n: n\ge N) \text{ where }\forall r :N> \deg f_r.
\]
The Calabi Yau property allows us to calculate the dimension of $iMj$:
\[
\Dim iMj = \Dim \Hom(P_i,Mj) = \Dim \Ext^3(S_i,Mj)\stackrel{CY}{=}\Dim \Hom(Mj,S_i) = \delta_{ij},  
\]
and conclude that $M$ must be isomorphic to the degree zero part of $A$. As $(2)$ implies there are only as many $f_r$ as there are arrows, we can conclude that the $f_r$ are linear and form a basis for $A_1$. Hence, by linearly combining our original relations, we can assume that the $f_r$ can be identified with the arrows. Let $r_a$ be the (nonzero) relation for which
$f_{r_a}=a$. This relation occurs only in the resolution of $S_{t(r_a)}=S_{h(a)}$ and therefore 
$h(a)=t(r_a)$ and $t(a)=h(r_a)$. 

Every arrow $a$ is contained in a cycle: $ar_a$, so if there is
a path between two vertices there is also a path in the opposite direction.
This means we that because $Q$ is assumed to be connected, $Q$ is also strongly connected.
We will now prove that all the $r_a$ have the same degree. 

Let $a$ be the arrow for which $r_a$ has minimal degree. 
First of all note that if two arrows $a,b$ share their heads then $\deg r_a = \deg r_b$ because they occur in the same resolution. Denote by $r_{ab}:=r_{a}b^{-1}$ the terms that appear in the middle map of the resolution. 
These terms are only nonzero if $t(b)=h(a)$. 
The fact that the maps in the resolutions form a complex implies that $\sum_{h(a)=i} ar_{ab}$ is zero in $A$. If $\deg r_a= \deg ar_{ab}$ is minimal then there exist scalars $(g_{bc})$ such that
\[
\sum_{h(a)=i} ar_{ab} = 
\sum_{\begin{smallmatrix}h(c)=h(b)\\t(c)=h(a)\end{smallmatrix}} g_{bc}r_c = \sum_{\begin{smallmatrix}h(c)=h(b)\\t(c)=t(b)\end{smallmatrix}} g_{bc}\sum_{t(d)=h(b)}r_{cd}d
\text{~~evaluated in $\C Q$.}
\] 
The $\deg r_c$ (which is the same for all $c$ with $h(c)=h(b)$ including $b$ itself) must also be minimal. All arrows following an arrow of minimal $r_a$-degree are also minimal, so by induction all arrows in $Q$ have the same degree.
 
We will now prove that $(g_{ab})$ can be seen as a diagonal matrix.
First note that
$$\Ext^1(S_i,S_j)=\Hom(\bigoplus_{t(a)=i} P_{h(a)},S_j) \cong \C^{\{i\to j\}}$$
and on the other hand
$$\Ext^2(S_j,S_i)=\Hom(\bigoplus_{t(r_a)=j} P_{h(r_a)},S_i)= \Hom(\bigoplus_{h(a)=j} P_{t(a)},S_i)\cong \C^{\{i\to j\}}.$$
We can compose the spaces in two different ways:
$$\Ext^1(S_i,S_j) \times \Ext^2(S_j,S_i)\to\Ext^3(S_j,S_j)\cong \C : (\xi_a)*(\eta_b) = \sum_a \xi_a\eta_a$$
and
$$\Ext^2(S_j,S_i) \times \Ext^1(S_i,S_j)\to\Ext^3(S_i,S_i)\cong \C : (\eta_b)*(\xi_a) = \sum_{a,b} g_{ab}\xi_a\eta_b$$

We only work out the last composition since the other one is similar.
We extend the sequence $(\eta_b)$ to a sequence running over all arrows by adding zeros.
We push out (dotted lines) the map $\eta$ forward along the resolution to obtain an exact sequence $S_i\to \dots\to S_j$:
\[
\xymatrix{
P_j \ar[r]^{\cdot c}&\bigoplus_{h(c)=j} P_{t(c)}\ar[r]^{\cdot r_{cd}}\ar[d]_{(\eta_d)}\ar@{:}[rd]&\bigoplus_{t(d)=j} P_{h(d)}\ar[r]^{\cdot d}\ar[d]\ar@{:}[rd]&P_j\ar[d]\ar@{:}[rd]\ar[r]&S_j\ar[d]\\
&S_i\ar[r]&\frac{S_i \bigoplus_{t(d)=j} P_{h(d)}}{((-\eta_c,r_{cd}), h(c)=j)}\ar[r]&
\frac{\bigoplus_{t(d)=j}P_{h(d)} \oplus P_j}{((-\delta_{cd},d), c: t(d)=j)}\ar[r]&S_j.}
\]
We use this sequence to pull back (dotted arrows) the map $(\xi_b)$
\[
\xymatrix@C=1.6cm{
P_i\ar[r]^{a}\ar@{.>}[rd]|m\ar@{.>}[d]|{\sum_{bc}g_{bc}\xi_b\eta_c}
&\bigoplus_{h(a)=i} P_{t(a)}\ar[r]^{r_{ab}}\ar@{.>}[rd]|{(0,\sum_br_{ab}\xi_b)}
\ar@{.>}[d]|{(0,\sum_b r_{ab}\xi_bd^{-1})}&\bigoplus_{t(a)=i} P_{h(b)}\ar[r]\ar@{.>}[rd]|{\xi_b}\ar@{.>}[d]|{(0,\xi_b)}&P_i\\
S_i\ar[r]&\frac{S_i \bigoplus_{t(d)=j} P_{h(d)}}{((-\eta_c,r_{cd}), h(c)=j)}\ar[r]&
\frac{\bigoplus_{t(d)=j}P_{h(d)} \oplus P_j}{((-\delta_{cd},d), c: t(d)=j)}\ar[r]&S_j
}
\] 
where 
\[
m = (0,\sum_{ab}ar_{ab}\xi_b d^{-1}) = (0,\sum_{bce}g_{bc}r_{ce}e\xi_b d^{-1})= (0,\sum_{bc}g_{bc}r_{cd}\xi_b)= (\sum_{bc}g_{bc}\eta_c\xi_b,0)
\]
Because of the Calabi Yau property there exist traces $\Tr_{S_j}: \Ext^{3}(S_j,S_j)\to \C$. As these Ext-spaces are one-dimensional we can represent these traces by scalars $\alpha_j$. Property \ref{ref-A.5.2-12} can be rewritten as
\se{
\Tr_{S_j}((\xi_a)*(\eta_b))&=\Tr_{S_i}((\eta_b)*(\xi_a))\\
\alpha_j  \sum_a \xi_a\eta_a &= \alpha_i \sum_{a,b} g_{ab} \xi_a\eta_b.
}
As this holds for arbitrary $(\xi_a)$ and $(\eta_b)$ we can conclude that
\[
g_{ab} = \frac {\alpha_{h(a)}}{\alpha_{t(a)}}\delta_{ab}.
\]

Now we construct the element
$$
\sum_{a,b \in Q_1} \alpha_{h(a)}ar_{ab}b,
$$ 
Which is a sum of cycles.
It is also a homogeneous element that is invariant under cyclic permutation:
\[
\sum_{a,b} \alpha_{h(a)}r_{ab}ab =\sum_{a,b} \alpha_{t(b)}r_{ab}ab= \sum_{a,b} \alpha_{t(b)}\frac{\alpha_{h(b)}}{\alpha_{t(b)}} br_{ba}a = \sum_{a,b} \alpha_{h(b)}br_{ba}a.
\]
This implies that we can identify it with $\cycle(W)$ where $W \in \C Q/[\C Q,\C Q]$ such that
$r_a$ is a scalar multiple of $\partial_a W$.

To prove the last condition on the structure of the quiver, assume first that $v$ is the tail of a unique arrow $a$ and let the $b_i$ be the vertices 
whose head is $t(a)$. As $r_a= \sum_i b_i r_{b_ia}$ and $r_a \ne 0$ in $\C Q$, there must be
at least one $r_{b_ia} \ne 0$ in $\C Q$ and because of its degree it is also nonzero in $A$.
Now $r_{b_ia}$ sits inside the kernel of $P_{h(a)} \stackrel{\cdot a}{\rightarrow} P_{t(a)}$ because $\partial_{b_i} W=r_{b_ia}a$. This would imply that the resolution for $S_{h(a)}$ is not exact.
Using right modules instead of left one proves that every vertex is also the tail of at least two arrows.
\end{proof}

For reasons of completeness we also include the descriptions of  
Calabi Yau algebras of smaller dimension because the techniques to do this are similar.

The zero-dimensional case is trivial and consists of the semi-simple algebras i.e. quivers without arrows. The one-dimensional case consists a direct sums of $\C[X]$ (disjunct unions of one-vertex-one-loop quivers). This is a consequence of property $C2$:
$\#\{i \leftarrow j \}=\Dim \Ext^1(S_i,S_j){=} \Dim \Hom(S_j,S_i) = \delta_{ij}$.

\begin{theorem}
If $A$ is Calabi Yau of dimension $2$ then $A$ is the preprojective algebra of a non-Dynkin quiver (for a definition of a preprojective algebra see \cite{CBHol}). 
\end{theorem}
\begin{proof}
As the global dimension of $A=\C Q/\cI$ is now $2$, the projective graded resolutions look like
\[
\xymatrix{
\bigoplus_{t(r)=i}P_{h(r)}\ar@{^(->}[r]^{\cdot ra^{-1}}&\bigoplus_{t(a)=i}P_{h(a)}\ar[r]^{\cdot a}&P_i\ar@{->>}[r]&S_i 
}
\]
From the Calabi Yau property C2, we deduce that
\[ 
\#\{r \in \cR|h(r)=i, t(r)=j\}=\Dim \Ext^2(S_i,S_j) \stackrel{CY}{=} \Dim \Hom(S_j,S_i) = \delta_{ij},
\]
i.e. for every vertex there is a unique relation and vice versa.

Now, similarly to the three dimensional case, we consider the finite dimensional quotient algebra $M = A/(ra^{-1}: r \in \cR, a\in Q_1, A_n: n\ge N)$  where $\forall r :N> \deg r$.
The Calabi Yau property allows us to calculate the dimension of $iMj$:
\[
\Dim iMj = \Dim \Hom(P_i,Mj) = \Dim \Ext^2(S_i,Mj)\stackrel{CY}{=}\Dim \Hom(Mj,S_i) = \delta_{ij}  
\]
and conclude that $M$ must be isomorphic to the degree zero part of $A$. This implies that the $ra^{-1}$ are all linear
and span $A_1$. For every $a$ there is also at most one $r$ such that $ra^{-1}$ is nonzero: the unique $r$ with
$t(r)=t(a)$. If we group the relations together in $R =\sum_{r \in \cR} r$ then   
there exists an invertible complex matrix $g_{ab}$ such that 
\[
Ra^{-1} = \sum_{a,b} g_{ab}b.
\]
We can use this $g$ to explicitly calculate the pairing (the calculation is analogous to the thee-dimensional case).
\[
\Ext^1(S_i,S_j) \times \Ext^1(S_j,S_i) \to \Ext^2(S_i,S_i) (\xi_a)*(\eta_b) = \sum_{ab}g_{ab}\xi_a\eta_b.
\]
Property $C3$ now implies that $g_{ab}$ is antisymmetric and non-degenerate
so using a base transformation on the arrows we can put $g_{ab}$ in its standard symplectic form.
The fact that $g_{ab}\ne 0 \implies h(a)=t(b) \wedge t(a)=h(b)$ indicates that this base transformation
only mixes arrows with identical head and tail. 
In this new basis the arrows can be partitioned in couples $(a,a^*)$ with $g_{aa^*} = 1$ and
$g_{ab}=0$ if $b\ne a^*$. The relation $R$ assumes the form of the standard preprojective relations:
\[
\sum_{a} aa^*-a^*a
\]
where $a$ runs over the unstarred half of the arrows.
Also $Q$ cannot be the double of a Dynkin quiver because $A$ must have global dimension $2$, see
\cite{CBHol}.
\end{proof}

\section{Selfdual Resolutions}

In this section we use the notion of selfdual resolutions to give a criterium to check whether a
vacualgebra $A_W$ is indeed Calabi Yau.

\subsection{Projective $A$-modules}
Let $A$ be a finitely generated graded algebra that is the quotient of a path algebra $\C Q$ and let $S=A_0$.

For every finite dimensional $S$-bimodule $T$ we can define a projective
$A$-bimodule
\[
F_T := A \otimes_{S} T \otimes_{S} A.
\]
We denote the full subcategory of $\Mod A-A$ containing these projective modules as $\cP$.
The basic objects of this category are of the form $F_{ij}:= F_{Si\otimes jS}=Ai \otimes jA$ with $i,j \in Q_0$.

The bimodule homomorphisms between $F_{T} \in \cP$ and a bimodule $M \in \Mod A-A$ can be identified with
\[
\Hom_{A-A}(F_{T},M)\cong T^* \stimes M.
\] 
The tensor product in this formula tensors over both the left and the right $S$ action.
The identification can be expressed explicitly as 
\[
\theta \stimes  m : b_1 \otimes_{S} t \otimes_{S} b_2 \mapsto \sum_{i,j \in  Q_0} \theta (itj) b_1 i  m j b_2       \]
A special role is played by $F_{S \otimes S}\cong A \otimes A$. We will denote this space by $F$. On this vector space we can define two commuting $A$-bimodule structures
\se{
F_{\text{Outer}}:&(a_1 (b_1 \otimes b_2)a_2)=a_1b_1\otimes b_2a_2,\\
F_{\text{Inner}}:&(a_1 (b_1 \otimes b_2)a_2)=b_1a_2\otimes a_1b_2.
}
If we use no subscript, we automatically assume the outer structure.
These structures are both isomorphic as bimodules to the free bimodule of rank one and the isomorphism between them is given by the twist
\[
\tau:F_{\text{Outer}} \to F_{\text{Inner}}:(b_1 \otimes b_2) \mapsto (b_2 \otimes b_1).
\]

The existence of these two commuting structures implies that for any $A$-bimodule $M$ the object $\Hom_{A-A}(M, F_{\text{Outer}})$ is again 
an $A$-bimodule using the inner structure. This bimodule will be denoted by $M^\vee$.
Maps can also be dualized in the standard way to turn $-^\vee$ into a functor:
\[
\forall f \in \Hom_{A-A}(M,N):\forall m \in M: \forall \nu \in N^\vee:
~f^\vee(\nu)(m) := \nu(f(m)). 
\]
For the standard projective bimodules we have the following natural identities
\begin{itemize}
\item 
$F_{T}^\vee = \Hom_{A-A}(F_{T},F)
\cong (T^* \stimes F)_{\text{Inner}}
\cong A \otimes_{S} T^* \otimes_{S} A
= F_{T^*},$
\item
$\Hom_{A-A}(F_{T},M) \cong T^* \stimes M \cong F_{T^*}\atimes M \cong F_{T}^\vee\atimes M.$
\end{itemize}
We can also write out the duality for the morphisms:
\[
(\theta \stimes a_1 \otimes_{S} t\otimes_{S} a_2)^\vee 
= t \stimes a_2 \otimes_{S} \theta \otimes_{S} a_1 
\]
These formulas imply that there is a natural equivalence between
\[
(- \atimes -)^* \text{ and }(-^\vee \atimes -^*): \cP \times \Rep A-A \to \Mod \C
\] 
and between $-^{\vee\vee}|_{\cP}$ and $\cP \hookrightarrow \Bimod A-A$. 
These functors and identities can be transferred to complexes if we assume that
\[
(M^\bullet)^* = (M_{-i}^*,-(d^{M}_{-i+1})^*)\text{ and }(P^\bullet)^\vee = (P_{-i}^\vee,-(d^{P}_{-i+1})^\vee)
\]
Keeping all this in mind we can propose the following definition:
\begin{definition}
A projective resolution $P^\bullet$ of left $A$-bimodules is selfdual with shift $n$ if and only if
there exists a commutative diagram
\[
\xymatrix{
P_n \ar[r]^{d_n}\ar[d]^{\alpha_n}&P_{n-1}\ar[d]^{\alpha_{n-1}}\ar[r]^{d_{n-1}}&\dots\ar[r]& P_1\ar[r]^{d_1}\ar[d]^{\alpha_1}&P_0\ar[d]^{\alpha_0}\\
P_0^\vee \ar[r]^{-d^\vee_1}&P_1^\vee\ar[r]^{-d^\vee_2}&\dots\ar[r]&P_{n-1}^\vee\ar[r]^{-d^\vee_n}&P_n^\vee
}
\]
for which the $\alpha_i$ are isomorphisms of $A$-bimodules. In short hand we can write $P^\bullet \cong (P^\bullet)^\vee[n]$.
\end{definition}

\begin{theorem}
If an algebra $A$ has a selfdual resolution of length $n$ with entries in $\cP$ then 
$A$ is Calabi Yau of dimension $n$.
\end{theorem}
\begin{proof}
Let $M^\bullet$ and $N^\bullet$ be two complexes in $\Rep A$. 
Standard homological algebra allows us to identify naturally
\se{
\Hom_{\cD^b\Rep A}(M^\bullet, N^\bullet)
&\cong \Hom_{\cD^b\Bimod A-A}(A,N^\bullet \otimes (M^\bullet)^*)\\
&\cong \Hom_{\cD^b\Bimod A-A}(P^\bullet,N^\bullet \otimes (M^\bullet)^*)\\
&\cong H^0\RHom_{\cD^b\Bimod A-A}(P^\bullet,N^\bullet \otimes (M^\bullet)^*)\\
&\cong H^0\Hom^\bullet_{\cK\Bimod A-A}(P^\bullet,N^\bullet \otimes (M^\bullet)^*)
}
So if we can prove that there is a natural equivalence between
\[
\Hom^\bullet_{\cK\Bimod A-A}(P^\bullet,N^\bullet \otimes (M^\bullet)^*)
\text{ and }
\Hom^\bullet_{\cK\Bimod A-A}(P^\bullet,M^\bullet \otimes (N^\bullet)^*[n])
\]
we are done.

Now using the fact that the resolution is composed of projectives in $\cP$ we can make the following identifications
\se{
&(\Hom^\bullet_{\cK\Bimod A-A}(P^\bullet,N^\bullet \otimes (M^\bullet)^*))^*\\
&\cong((P^\bullet)^\vee \atimes N^\bullet \otimes (M^\bullet)^*)^*\\
&\cong (P^\bullet)^{\vee\vee} \atimes (N^\bullet \otimes (M^\bullet)^*)^*\\
&\cong (P^\bullet) \atimes M^\bullet \otimes (N^\bullet)^*\\
&\stackrel{\alpha}{\cong}  (P^\bullet)^\vee[n] \atimes M^\bullet \otimes (N^\bullet)^*\\
&\cong  (P^\bullet)^\vee \atimes M^\bullet \otimes (N^\bullet)^*[n]\\
&\cong \Hom^\bullet_{\cD^b\Bimod A-A}(P^\bullet,M^\bullet \otimes (N^\bullet)^*[n])
}
which are natural in the $M^\bullet$ an $N^\bullet$.
\end{proof}

For an explicit write-out of the corresponding pairing between $\Hom_{\cK\Bimod A-A}(P^\bullet,N^\bullet \otimes (M^\bullet)^*)$ and 
$\Hom_{\cK\Bimod A-A}(P^\bullet,M^\bullet \otimes (N^\bullet)^*[n])$ we first need some notation: for simplicity we will work with
elements that are pure tensors:
\se{
f \in \Hom_{\cK\Bimod A-A}(P^\bullet,N^\bullet \otimes (M^\bullet)^*): f^{ij} &= \phi^{ij} \atimes \mu^{ij} \in {P^i}^\vee \atimes N^{i+j}\otimes {M^{j}}^*\\
g \in \Hom_{\cK\Bimod A-A}(P^\bullet,M^\bullet \otimes (N^\bullet)^*[n]): g^{ij} &= \gamma^{ij} \atimes m^{ij} \in {P^i}^\vee \atimes M^{j} \otimes (N^{n-i+j})^*.
}
With these expressions for $f$ and $g$ we can track back the pairing in the previous identifications:
\[
\<f,g\>_{M^\bullet N^\bullet}=   \sum_{ij} \Tr \mu^{ij}\circ \phi^{ij}\alpha_{n-i}^{-1}(\gamma^{n-i,j})m^{n-i,j}.
\]

\subsection{Superpotentials and Selfduality}

In the case of a graded algebra $A := \C Q / \cI,~\cI \subset \cJ^2$  one can construct its minimal resolution 
using standard presentations of $\cI^n/\cI^{n+1}$. These objects, introduced in \cite{kingbutler}, consist of quintuples
$(U,V,r,l,\Delta)_n$ where
\begin{enumerate}
\item
$U,V \subset \cI^n$ are $S$-bimodule complements such that
\[\cI^n = U \oplus \cJ\cI^n+\cI^n\cJ \text{ and }\cJ\cI^n\cap \cI^n\cJ = V \oplus \cJ\cI^n\cJ,\]
\item
$r,l: \cI^n \to A \otimes_S U \otimes_S A$ are a $\C Q-S$ and a $S-\C Q$-bimodule section of the $\C Q$-bimodule morphism
\[
e : A \otimes_S U \otimes_S A \to \cI^n: 1 \otimes_S u \otimes_S 1 \mapsto u.
\]  
and use these to define a map
\[
d : A \otimes_S V \otimes_S A  \to A \otimes_S U \otimes_S A: 1 \otimes_S v \otimes_S 1 \mapsto l(v) - r(v)
\]
\item
$\Delta: \cI^{n} \to A \otimes_S V \otimes_S A$ is a $\C Q$-bimodule derivation (i.e. a $S$-bimodule morphism satisfying $\Delta(azb) = \Delta(az)b +a\Delta(zb)- a\Delta(z) b$)  such that
$d\Delta=l - r$ and $\forall v \in V: \Delta(v) = 1 \otimes_S v \otimes_S 1$.
\end{enumerate}
Although the map $d$ is a morphism as $\C Q$-bimodules it can also be considered as a morphism of $A$-modules $d_A: F_V \to F_U$ because the $\C Q$-action factors over $A$. The same can be done with $e$ provided we factor out $\cI^{n+1}$ in the target: $e_A: F_U \to \cI^n/\cI^{n+1}$.
To turn $\Delta$ into a $A$-bimodule morphism we have to do two things: look at the subspace $\cI^{n+1} \subset \cI^{n}$ (this turns the derivation law into a morphism law) and mod out $\cI^{n+2}$ (this turns the domain into a $A$-bimodule):
\[
c_A : \frac{\cI^{n+1}}{\cI^{n+2}} \to F_V : x + \cI^{n+2} \mapsto \Delta(x).
\] 
These maps can be packed together in sequences of $A-A$ bimodules
\[
\xymatrix
{
0\ar[r]&\frac{\cI^{n+1}}{\cI^{n+2}}\ar[r]^{c_A}&F_{\frac{\cJ\cI^n\cap \cI^n\cJ}{\cJ\cI^n\cJ}}\ar[r]^{d_A}
&F_{\frac{\cI^n}{\cJ\cI^n+\cI^n\cJ}}\ar[r]^{e_A}&\frac{\cI^{n}}{\cI^{n+1}}\ar[r]&0.
}
\]
In \cite{kingbutler} it is proved that these sequences are exact and they can be spliced together
to get a projective bimodule resolution of $\cI^0/\cI^1=A$. This resolution is not minimal but it can
be made minimal if one cuts out the excess summands that occur at the splicing boundaries. These terms are of the form
\[
A \otimes_S \frac{\cI^{n+1}}{\cI^{n+1}\cap \cJ \cI^{n} \cJ} \otimes_S A \cong
A \otimes_S \frac{\cI^{n+1}+\cJ \cI^{n} \cJ}{\cJ \cI^{n} \cJ} \otimes_S A \subset F_{\frac{\cJ\cI^n\cap \cI^n\cJ}{\cJ\cI^n\cJ}}
\]

We will now apply this to the case of Calabi Yau algebras of dimension $3$.
As we already know from section \ref{nodige} the ideal is generated by
\[
\partial_a W, a \in Q_1
\]
where $W \in \C Q/[\C Q,\C Q]$ is a superpotential and as the global dimension is $3$ we only need to 
look at the standard presentations for $n=0,1$.

The case $n=0$ has the same form for every algebra
\begin{itemize}
\item
$U_0=S$, $V_0=\C Q_1$
\item
$l_0 : a \mapsto a \otimes_S t(a) \otimes_S 1$\\
$r_0 : a \mapsto 1 \otimes_S h(a) \otimes_S a$,
\item
$\Delta : a_1\cdots a_k \mapsto \sum_{1\le j \le k} a_1\cdots a_{j-1}\otimes_S a_j \otimes_S a_{j+1}\cdots a_k).$
\end{itemize}
For $n=1$ we do not need to bother about the $\Delta_1$ because it does not affect the minimal resolution:
\begin{itemize}
\item
$\cI = \C\{\partial_a W, a \in Q_1\} \oplus \cJ\cI+\cI\cJ,~ U_1= \C\{\partial_a W, a \in Q_1\}\cong \C Q^{op}_1$.\\
For $V_1$ we chose a complement that contains the subspace $\C\{i\cycle W, i \in Q_0\}$.
\item
$l_1 : x\partial_a Wy \mapsto x \otimes_S \partial_a W \otimes_S y$ if $y \notin \cI$,\\
$r_1 : x\partial_a Wy \mapsto x \otimes_S \partial_a W \otimes_S y$ if $x \notin \cI$.
\end{itemize}
Because the $i\cycl W$ are not contained in $\cI^2$ they are not cut out by restricting to the
minimal resolution. Moreover, because 
\se{
\Ext^3_A(S,S)&=\Hom_S(\frac{V_1}{\cI^2\cap V_1},S)\\
&\stackrel{CY}{\cong} \Hom_A(S,S)^* \\
&\cong \Hom_S(\C\{i, i\in Q_0\},S)^*\\ 
&\cong \Hom_S(\C\{i\cycl W, i\in Q_0\},S)^* 
}
We have that the third term in the minimal resolution must be $F_{\frac{V_1}{\cI^2\cap V_1}}=F_{\C\{i\cycl W, i\in Q_0\}}\cong F_S$.
Putting everything together we get 
\[
\xymatrix
{
F_{\C\{i\cycl W, i\in Q_0\}}\ar[r]^{\delta_3}&
F_{\C\{\partial_a W, a\in Q_1\}} \ar[r]^{\delta_2}&
F_{\C Q_1} \ar[r]^{\delta_1}&
F_{S} \ar[r]&0
}
\]
with maps
\se{
\delta_1(1 \otimes_S a \otimes_S 1) &= a \otimes_S t(a) \otimes_S 1  -  1 \otimes_S h(a) \otimes_S a\\
\delta_2(1 \otimes_S \partial_a W \otimes_S 1) &= \Delta(\partial_a W)\\
\delta_3(1 \otimes_S W \otimes_S 1) &= \sum_{a \in Q_1} a \otimes_S \partial_a W \otimes_S 1 - 1 \otimes_S \partial_a W \otimes_S a
}

A more explicit write-out of the complex $C_W$, whose $0^{th}$ homology is equal to $A$, in terms of
the basic projective $F_{ij}$
looks like
\[
C_W:\hspace{.5cm}\xymatrix{
\underset{i \in Q_0}\bigoplus F_{ii}\ar[r]^{(\cdot \tau d a \cdot)}&
\underset{a \in Q_1}\bigoplus F_{t(a)h(a)}\ar[r]^{(\cdot\partial^2_{ba}W \cdot)}&
\underset{b \in Q_1}\bigoplus F_{h(b)t(b)}\ar[r]^{(\cdot db \cdot)}&
\underset{i \in Q_0}\bigoplus F_{ii}\ar[r]^m&A}
\]
where the differential is $da := a\otimes t(a) - h(a) \otimes a$ and the second derivatives 
are $\partial^2_{ba}W= \pi_{F_{t(b)h(b)}}\Delta \partial_a W$. More explicitly, if $c$ is a cycle then
\[
\partial^2_{ba} c = \sum_{p_1,p_2: \cycle ap_1bp_2=\cycl c} p_1 \otimes p_2.
\]
Note that because $\cycl W$ is invariant under cyclic permutation, $\partial^2_{ab}W = \tau\partial^2_{ba}W$  

\[
\xymatrix@R=.3cm{
\left(\underset{i \in Q_0}\bigoplus F_{ii}\right.\ar[r]^{(\cdot \tau d a \cdot)}&
\underset{a \in Q_1}\bigoplus F_{t(a)h(a)}\ar[r]^{(\cdot\partial^2_{ba}W \cdot)}&
\underset{b \in Q_1}\bigoplus F_{h(b)t(b)}\ar[r]^{(\cdot db \cdot)}&
\left.\underset{j \in Q_0}\bigoplus F_{jj}\right)^\vee\\
=\underset{j \in Q_0}\bigoplus F_{jj}^\vee\ar[r]^{(\cdot d b \cdot)^\vee}\ar[d]^\tau&
\underset{b \in Q_1}\bigoplus F_{h(b)t(b)}^\vee\ar[r]^{(\cdot\partial^2_{ba}W \cdot)^\vee}\ar[d]^\tau&
\underset{a \in Q_1}\bigoplus F_{t(a)h(a)}^\vee\ar[r]^{(\cdot \tau da \cdot)^\vee}\ar[d]^\tau&
\underset{i \in Q_0}\bigoplus F_{ii}^\vee\ar[d]^\tau\\
=\underset{j \in Q_0}\bigoplus F_{jj}\ar[r]^{(\cdot \tau d b \cdot)}&
\underset{b \in Q_1}\bigoplus F_{t(b)h(b)}\ar[r]^{(\cdot\partial^2_{ab}W \cdot)}&
\underset{a \in Q_1}\bigoplus F_{h(a)t(a)}\ar[r]^{(\cdot da \cdot)}&
\underset{i \in Q_0}\bigoplus F_{ii}}
\]
This complex is selfdual and the isomorphism connecting the complex with its dual is composed of
the standard identifications we used in the previous paragraph.

So the sufficient condition of selfduality is also necessary for Calabi Yau algebras of dimension $2$.
\begin{theorem}
A vacualgebra $A_W$ is Calabi Yau of dimension $3$ if and only if the complex $C_W$ is a projective resolution of $A_W$ as an $A_W$-bimodule. 
\end{theorem}

This fact has a nice interpretation for the classification of \emph{good superpotentials} i.e. superpotentials with a vacualgebra that is indeed Calabi Yau.
\begin{corollary}
For a given quiver $Q$ and a given dimension $d$, the subset of $\Sup_d Q$ of good superpotentials of degree $d$ is either the empty set or almost everything (in the measure theoretic sense). 
\end{corollary}
\begin{proof}
The condition we must check is that the standard complex is indeed a resolution. Because the resolution is graded we can check this separately for every degree so the subspace of good superpotentials is 
an intersection of a countable number of Zariski open sets. If one of these sets is empty we're in the first case and otherwise the complement of this set is a countable union of hypersurfaces, which has measure zero for the standard measure on $\C^n$.
\end{proof}

\begin{remark}
For global dimension two we can do a similar thing. Recall that if $A$ is Calabi Yau of dimension two, then
the set of arrows partitions in pairs $(a,a^*)$ with opposite head and tail.  

The selfdual resolution now looks like
\[
\hspace{.5cm}\xymatrix@C=1.5cm{
\underset{i \in Q_0}\bigoplus F_{ii}\ar[r]^<<<<<<<<{\left(\begin{smallmatrix}\cdot \tau da^* \cdot\\ \cdot -\tau da \cdot\end{smallmatrix}\right)}&
\underset{(a,a^*)}\bigoplus F_{t(a)h(a)}\oplus F_{t(a^*)h(a^*)} \ar[r]^>>>>>>>>{\left(\begin{smallmatrix}\cdot da \cdot\\\cdot da^* \cdot\end{smallmatrix}\right)}&
\underset{i \in Q_0}\bigoplus F_{ii}\ar[r]^m&A}
\]
This is indeed
the standard resolution for preprojective algebras of non-Dynkin quivers (see \cite{CBHol}).
\end{remark}
\subsection{The matrix valued Hilbert polynomial}
For a graded algebra $A=\C Q/(\cR)$ one can define the \emph{matrix valued Hilbert series}
\[
H_A(t) := h_0 + h_1 t + h_2t^2 + \cdots
\]
where the $h_k$ are matrices in $\Mat_{\#Q_0 \times \#Q_0}(\C)$ and
\[
(h_k)_{ij} = \dim iA_kj~\text{($A_k$ is the degree $k$ part of $A$)}
\]

The matrix valued Hilbert series of a Calabi Yau algebra can be computed from its bimodule resolution: 
\begin{theorem}
If a vacualgebra $A_W$ with $\deg W= d \ge3$ is Calabi Yau then  
\[
H_{A_W}(t) = \frac{1}{1 - M_Q t + M_Q^Tt^{d-1} - t^d}
\]
where $M_Q$ is the incidence matrix of $Q$.
This equality must be evaluated in the ring of formal power series
$\Mat_{\#Q_0 \times \#Q_0}(\C[\![t]\!])$.
\end{theorem}
\begin{proof}
The Hilbert polynomial of $F_{kl}$ is equal to 
\[
H_{F_{kl}}(t)= H_A(t) e_{kl} H_A(t)
\]
where $e_{kl}$ is the matrix with $1$ on the entry $k,l$ and zero elsewhere.
So from the exactness of the resolution $C_W$ and the fact that $H_0(P^\bullet)=A$ we get
\se{
H_A(t) &= H_{F_S} - t(H_{F_{\C Q_1}} - t^{d-2}(H_{F_{\C Q_1^{op}}} - tH_{F_{S}}))\\
&= H_A(t) 1 H_A(t) - t H_A(t) M_Q H_A(t) + t^{d-1}H_A(t) M_Q^{T} H_A(t) - t^d H_A(t) 1 H_A(t)
}
Note that $H_A(t)$ is invertible because $H_A(0)=\id$.
Multiplying to the left and the right by $H_A(t)^{-1}$ and taking the inverse we obtain
the equality.
\end{proof}

The bimodule resolution gives us also resolutions of the left modules $S$.
Writing out the dimensions of these resolutions
gives the equation
\[
1 =  H_A(t) - t  M_Q H_A(t) + t^{d-1}  M_Q^{T} H_A(t) - t^d H_A(t).
\]
This is nothing new, but as this equation corresponds to a real resolution we can derive certain inequalities 
that must be met:
\begin{itemize}
\item[I1]
$H_A(t) \ge 0$
\item[I2]
$(M_Q^{T} - t) H_A(t) \ge 0$
\item[I3]
$(M_Q - M_Q^{T}t^{d-2} - t^{d-1}) H_A(t) \ge 0$
\end{itemize}
Note that a matrix valued series $f(t)$ is positive if all its entries $(f_k)_{ij}$ are positive.
These inequalities can be useful to check whether quivers have good superpotentials of a given degree.

\section{Applications}
\subsection{Groebner Bases and Superpotentials}
To show that for a given quiver and a given degree there exist good superpotentials one has to check 
whether one can find a superpotential $W$ such that $C_W$ is exact. 
To do this we will use the technique of Groebner bases as outlined in \cite{groeb}, adapted to path algebras.
Suppose $Q$ is a quiver with $n$ arrows and put an order on the arrows: $a_1>\cdots> a_n$. One can extend this
order to the set of paths with nonzero length 
using the \emph{deglex ordering} method:
\[
a_{i_1}\cdots a_{i_p} < a_{j_1}\cdots a_{j_q}
\]
if and only if $p<q$ or $p=q$ and $\exists \nu \le p:a_{i_\nu}< a_{j_\nu} \wedge \forall \mu < \nu: i_\mu=j_\mu$.
We denote the leading monomial term (according to the deglex ordering) of $f \in \C Q$ by $\lt(f)$. Recall
that a (not necessarily finite) set of elements $G \subset \cI\lhd \C Q$ is a \emph{Groebner basis} if all $\lt(g), g\in G$ are different and 
\[
\lt(\cI) := (\lt(f): f\in \cI) = (\lt(g): g \in G)  
\] 
where the equality is taken as ideals in $\C Q$.
Groebner bases are very useful in determining the structure of an algebra.
They can be used to determine the Hilbert polynomial because
\[
H_{\C Q/\cI} = H_{\C Q/\lt\cI}
\]
and they can be used to check whether certain expressions in $\C Q$ are 
zero in $\C Q/\cI$:
\[
f \in \cI \implies \lt(f) \in \lt(\cI) =(\lt(g): g \in G)
\]
To check whether a given set of relations is indeed a Groebner basis 
one can use the method of Bergman diamonds \cite{bergman}.
For any $f$ in $\C Q$, an elementary reduction of $f$ by $g \in G$ is the new expression
\[
\rho_{g}(f) := f - \zeta agb \text{ if $a,b$ are paths and $\zeta \in \C$ s.t.}\lt(f)=\zeta a\lt(g)b\text{ or $f$ otherwise}.
\]
If $G$ is a set of relations then a triple of monomial terms $(a,b,c)$ is called
an \emph{ambiguity} if $ab=\lt(g_1), bc=\lt(g_2)$ with $g_1,g_2 \in G$. An ambiguity is called resolvable if there is a sequence of elementary reductions such that
\[
\rho_1\cdots\rho_m(g_1c -ag_2)=0.
\]
Now Bergman's diamond lemma states that if all leading terms are different and 
all ambiguities are resolvable then $G$ is a Groebner basis. 

We will now give a useful criterion to find good superpotentials.
\begin{lemma}\label{condsup}
Suppose every vertex in $Q$ is the source and the target of at least two arrows and
$W$ is a superpotential such that
\begin{itemize}
\item
The leading terms of the relations $\partial_a W$ are all different and the ambiguities are in $1$ to $1$ correspondence to the vertices and are of the form
\[
\mathrm{amb_v}=(a, \lt(\partial_a W)b^{-1}, b )=(a, a^{-1} \lt(\partial_b W), b )\text{ with } a\lt(\partial_a W)=\lt(\partial_b W)b=\lt(vWv),
\]
\item
for every vertex $v$ there is at least one arrow $a, t(a)=v$ such that $\forall b \in Q_0:\lt(\partial_b W)a^{-1}=0$.
\end{itemize}
then $A_W$ is Calabi Yau.
\end{lemma}
\begin{proof}
First note that the condition implies that $\{\partial_a W: a \in Q_1\}$ is a Groebner basis:
an ambiguity of the form $\mathrm{amb}_v(a, \lt(\partial_a W)b^{-1}, b)$ is resolvable
because 
\[
\sum_{h(c)=v} c \partial_{c} W = \sum_{t(c)=v} \partial_{c} W c 
\]
and hence
\[
a\partial_{a} W - \partial_{b}W b =  \sum_{t(c)=v,c \ne b} \partial_{c} W c- \sum_{h(c)=v, c \ne a } c \partial_{c} W.  
\]
Note that the leading terms of the summands in the right hand side are all different because the $\lt(\partial_c W)$ are and there is only one ambiguity corresponding to $v$. We can remove each term using an elementary reduction, starting with the one with the highest
leading term. Therefore the ambiguity is resolvable.

To calculate the Hilbert series one must calculate $(h_k)_{vw}$, which is equal to the number of words between $v$ and $w$ of a given length $k$ not containing $\lt(\partial_a W)$'s. This can be done using recursion:
\[
(h_k)_{vw} = \underbrace{\sum_u h^{k-1}_{vu}\#\{u\leftarrow w\}}_{\text{add an arrow}} - 
\underbrace{\sum_u h^{k-d+1}_{vu}\#\{w\leftarrow u\}}_{\text{remove those ending in $\lt(\partial_b W)$}} +
\underbrace{h^{k-d}_{vw}}_{\text{remove double counting}}.
\]
There are no further terms needed: a word ending $w$ can only be double counted once because of the form
and number of the ambiguities. The Hilbert series of $A_W$ is thus
\[
H_{A_W}(t) = \frac{1}{1 - M_Qt + M_Q^Tt^{d-1} -t^d}.
\]
Using the exactness of the first 2 terms of $C_W$ we can calculate the Hilbert series of the kernel of the third map
\[
H_{A_W}(t) - kt H_{A_W}(t)^2 + kt^{d-1} H_{A_W}(t)^2 = t^d H_{A_W}(t)^2.
\]
This is the same as the Hilbert series of the last term so if we can prove that the last map is an injection we are done.
There is indeed no element $\sum_j f_j\otimes g_j \in A_W\otimes_S A_W$ such that
\[
\forall i \ge k: \sum_j f_jb\otimes g_j - f_j\otimes bg_j  = 0.
\] 
The deglex ordering on $\C Q$ can be transferred to and ordering on the monomials of  $\C Q\otimes_S \C Q$:
\[
v_1 \otimes v_2 > w_1 \otimes w_2 \iff v_1>w_1 \text{ or }v_1=w_1\text{ and }v_2>w_2.  
\]
This ordering is compatible with the multiplicative structure on $\C Q\otimes_S \C Q$.
Let $f_1 \otimes g_1$ be the highest order term, then the highest order term
of $\sum_j f_jb\otimes g_j - f_j\otimes bg_j$ is $f_1b\otimes g_1$. 
Therefore $f_1 \notin (\lt(\partial_{a}W):a \in Q_1)$ but $f_1b \in (\lt(\partial_{a}W):a \in Q_0)$ for every $b \in Q_0$.
This would imply that for every $b$ with $h(b)=t(f_1)$ there is a $c \in Q_1$ such that $f_1b$ ends in 
$\lt(\partial_{c}W)$, contradicting the second condition on $W$.
\end{proof}
The conditions imposed on the superpotential are very strict and there are far more good superpotential that do not
meet these conditions. In general the ideal generated by a good superpotential will not have a finite Groebner basis.
However for many quivers and degrees we will be able to find superpotentials that satisfy the demands of the lemma.

\subsection{The one vertex situation}
First note that if $Q$ has only one vertex and one loop, then none of the vacualgebras can be Calabi Yau of dimension $3$
because these algebras are finite dimensional and hence $H_A(t)$ cannot be the inverse of the polynomial $1 - t+ t^{d-1} - t^d$.

So, in this section, let $Q$ be a quiver with one vertex and $k\ge 2$ loops and let $\Sup_d \subset \C Q/ [\C Q,\C Q]$ be the subspace of all superpotentials of degree $d$ with $d\ge 3$. We will show that the space of good superpotentials is non-empty if and only if $(k,d) \ne (2,3)$.

If $(k,d)=(2,3)$ then there are no good superpotentials because the inequality $(I2)$ does not hold:
\[
(2-t)\frac{1}{1 - 2t+2t^2 -t^d} = 2 + 3 t + 2 t^2  - t^3  + 2 t^4  + \cdots \not \ge 0. 
\]

For every other couple $(k,d)$ we can find at least one good superpotential.

\begin{lemma}
Take $\C Q \cong \C\<X_1,\dots,X_n\>$ and $X_1>X_2>\cdots > X_n$, then the following superpotentials are good:
\begin{enumerate}
\item $W = X_1X_2X_3 + X_1X_3X_2 + \sum_{j>3}X_1X_j^2 +[\C Q,\C Q]$,
\item $W = \sum_{k\ge l > 1}X_1^{d-2}X_lX_k +[\C Q,\C Q]$.
\end{enumerate}
\end{lemma}
\begin{proof}
We calculate the leading terms of the relations
\begin{enumerate}
\item
$\lt(\partial_{X_1} W) =X_2X_3,~ \lt(\partial_{X_2} W)=X_1X_3,~ \lt(\partial_{X_3} W)X_1X_2,~\lt(\partial_{X_k} W) X_1X_k,\dots$ 
\item
$\lt(\partial_{X_1} W) = X_1^{d-3}X_2^2,~\lt(\partial_{X_2} W) = X_1^{d-1}X_2,~\dots,~\lt(\partial_{X_k} W) = X_1^{d-1}X_k$.
\end{enumerate}
The only ambiguity we can construct is 
\begin{enumerate}
\item  $(X_1,X_2,X_3)$ between $\partial_{X_1} W$ and $\partial_{X_2} W$,
\item $(X_1,X_1^{d-3}X_2,X_2)$ between 
$\partial_{X_1} W$ and $\partial_{X_2} W$.
\end{enumerate}
We also see that none of the leading terms ends in $X_1$.
\end{proof}
\begin{remark}
In the cases where $(k,d)$ equals $(2,4)$ or $(3,3)$ one can obtain a complete classification
of the good superpotentials because then we are in the case of Artin-Shelter regular algebras \cite{AS}.
\end{remark}
\subsection{Special Quivers}
The simplest quivers with more than one vertex that can have good potentials are
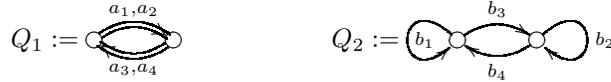
\[
Q_1 := \xymatrix{\vtx{}\ar@2@/^/[r]^{a_1,a_2}&\vtx{}\ar@2@/^/[l]^{a_3,a_4}}
\hspace{2cm}
Q_2 := \xymatrix{\vtx{}\ar@/^/[r]^{b_3}\ar@(lu,ld)^{b_1}&\vtx{}\ar@/^/[l]^{b_4}\ar@(ru,rd)^{b_2}}
\]
\begin{theorem}
\begin{itemize}
\item[]
\item
$\C Q_1/[\C Q_1,\C Q_1]_d$ contains good superpotentials if and only if $d\ge 4$ and $d$ is even.
\item
$\C Q_1/[\C Q_2,\C Q_2]_d$ contains good superpotentials if and only if $d\ge 4$. 
\end{itemize}
\end{theorem}
\begin{proof}
For both quivers $d$ must be bigger than or equal to $4$ because otherwise the inequalities $I1-I3$ are not satisfied. For $Q_1$, $d$ must be even because every cycle has even length.

Assume the orders $a_1>a_2>a_3>a_4$, $b_1>b_2>b_3>b_4$ and define the following superpotentials
\se{
Q_1&: a_1a_3(a_2a_4)^{\frac d2 -1} + a_3a_1(a_4a_2)^{\frac d2 -1}+[\C Q,\C Q]\\
Q_2&: b_1^{d-2}b_3b_4+b_2^{d-2}b_4b_3+[\C Q,\C Q]
}
The leading terms of the relations are
\se{
Q_1&:a_3(a_2a_4)^{\frac d2-1},~
a_3a_1a_4(a_2a_4)^{\frac d2-2},~
a_1(a_4a_2)^{\frac d2-1},~
a_1a_3a_2(a_4a_2)^{\frac d2-2},~\\
Q_2&:
b_1^{d-3}b_3b_4,~
b_2^{d-3}b_4b_3,~
b_2^{d-2}b_4,~
b_1^{d-2}b_3
}
For each of the quivers there are two ambiguities (one for each vertex)
\se{
Q_1&:(a_1,a_3a_2(a_4a_2)^{\frac d2-2},a_4),~ (a_3,a_1a_4(a_2a_4)^{\frac d2-2},a_2)\\
Q_2&:(b_1,b_1^{d-3}b_3,b_4),~(b_2,b_2^{d-3}b_4,b_3)
}
Finally none of the relations end in $a_1,a_3$ and $b_1,b_2$.
\end{proof}
The method described above can be extended to lots of other quivers and degrees, especially quivers
of the form 
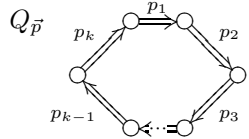
\[
Q_{\vec p}\xymatrix@=.5cm{
&\vtx{}\ar@2[r]^{p_1}& \vtx{}\ar@2[rd]^{p_2}&\\
\vtx{}\ar@2[ru]^{p_k}&&&\vtx{}\ar@2[ld]^{p_3}\\
&\vtx{}\ar@2[lu]^{p_{k-1}}& \vtx{}\ar@2[l]|\cdots&
}
\]
The number of arrows between consecutive vertices can differ (but is $\ge 2$).
\begin{theorem}
Let $Q$ be a quiver of the form above with $k\ge 2$ vertices and let $p_i\ge 2$ be the number of
arrows between the $i^{th}$ and the ${i+1}^{th}$ vertex. 
If $d=\ell k$ with $\ell \ge 2$ then 
$Q$ has good superpotentials of dimension $d$.
\end{theorem}
\begin{proof}
For every vertex $v \in Q_0$, we will denote the consecutive vertex by $v+1$, so $\forall a \in Q_1: h(a)=t(a)+1$.
Fix an order on the arrows of $Q$ and let $a_i,b_i$ the highest and second highest arrow arriving in the vertex $i$.

Define the superpotential
\se{
W := &\sum_{i \in Q_0} a_ia_{i-1}b_{i-2}\cdots b_{i-\ell k+1}\\
+ &\sum_{c \ne a_i,b_i} ca_{h(c)-1}b_{h(c)-2}\cdots b_{h(c)-k +1}(cb_{h(c)-1}b_{h(c)-2}\cdots b_{h(c)-k +1})^{\ell - 1}+[\C Q,\C Q]
}
The leading terms of the relations now look like
\se{
\lt(\partial_{a_i} W) &= a_{i-1}b_{i-2}\cdots b_{i-\ell k+1}\\
\lt(\partial_{b_i} W) &= a_{i-1}a_{i-2}b_{i-3}\cdots b_{i-\ell k+1}\\
\lt(\partial_{c} W) &= a_{h(c)-1}b_{h(c)-2}\cdots b_{h(c)-k +1}(cb_{h(c)-1}b_{h(c)-2}\cdots b_{h(c)-k +1})^{\ell - 1}
}
It is easy to check that all ambiguities are of the form $(a_i,a_{i-1}b_{i-2}\cdots b_{i-\ell k+2},b_{i-\ell k+1})$  
and none of the relations ends in some $a_i$.
\end{proof}
\begin{remark}
If $\ell=1$ the situation is more complicated because the solutions of the inequalities I1-I3 are not easy to determine. 
It is not the case that if they are satisfied for $Q_{\vec p}$ that they are also satisfied for a quiver $Q_{{\vec p}'}$
with $(p_1',\cdots,p_k')\ge (p_1,\cdots, p_k)$. F.i. a quiver with arrows $\vec p = (2,2,2,2)$ has good superpotentials
but one with $\vec p =(6,2,2,2)$ has not.
\end{remark}

The method of finding these very special superpotentials does not always work. As a counterexample consider the quiver
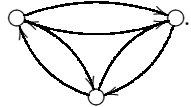
\[
\xymatrix{
\vtx{}\ar@/^/[rr]\ar@/^/[rd]&&\vtx{}\ar@/^/[ll]\ar@/^/[ld]\\
&\vtx{}\ar@/^/[ru]\ar@/^/[lu]&
}.
\]
One can check that there are no superpotentials of dimension $4$ 
satisfying the conditions from \ref{condsup} although Groebner basis computations in GAP \cite{gap}(up to a certain degree because the full Groebner basis could be infinite) seem to indicate that a generic superpotential is indeed good.

The general picture that arises from computations is that as soon as conditions I1-I3 are met by the Hilbert series then there do exist good superpotentials, but we have no proof for this statement.

\section{Aknowledgements}
This paper arose from discussions I had during an informal seminar in the summer of 2005. I would like to thank the other regular participants, Geert van de Weyer, Koen De Naeghel, Adam-Christiaan Van Roosmalen, Joost Vercruysse, Tor Lowen and Stijn Symens for their helpful
comments. I would especially like to thank Michel Van den Bergh for sharing his useful insights and for providing the appendix.

\newpage
\appendix
\section{The signs of Serre duality}
\section*{by Michel Van den Bergh}
\subsection{Introduction}
In this self-contained appendix we determine the exact signs which
occur in Serre duality (see for example Proposition \ref{ref-A.5.2-12}
for the Calabi-Yau case).  Although the answer is the obvious, the
verification turned out to be slightly more tricky than foreseen.

We thank Bernhard Keller for pointing out Example \ref{ref-A.3.2-2} (see \cite{Keller} and \cite{Verdier} for further information)
and suggesting that,
likewise, the correct
signs in Serre duality should be determined by the requirement that the Serre
functor be exact.
\subsection{Graded categories}
\begin{definitions} A graded (pre-additive) category is a pair $(\Cscr,S)$ where
$\Cscr$ is a pre-additive category and $S$ is an automorphism of $\Cscr$.
\end{definitions}
\begin{remarks} It is customary to only require $S$ to be an
  autoequivalence. The stronger condition that $S$ is an automorphism is
  usually satisfied in practice and up to an appropriate notion of
  equivalence we may always reduce to this case.
\end{remarks}

In a graded category $(\Cscr,S)$ we may define the \emph{graded $\Hom$-sets}
between objects by
\[
\Hom^i_\Cscr(A,B)=\Hom_\Cscr(A,S^i B)
\]
and 
\[
\Hom^{\text{gr}}_\Cscr(A,B)=\bigoplus_i \Hom^i_\Cscr(A,B)
\]
There is an obvious graded composition
\se{
\label{ref-1-0}
-\ast-:\Hom^j_\Cscr(B,C)\times \Hom_\Cscr^i(A, B)\r
\Hom^{i+j}_\Cscr(A,C):(g,f)\mapsto:
S^i(g)f
}
We denote by $\Cscr^{\text{gr}}$ the category $\Cscr$ equipped with graded
$\Hom$-sets.

\medskip

A \emph{graded functor} between graded categories $(\Cscr,S)$,
$(\Dscr,T)$ is an additive  functor $U:\Cscr\r \Dscr$ together with 
a natural isomorphism $\eta^U:U\circ S\r T\circ U$. By a slight
abuse of notation we will write the composition
\[
U\circ S^i\r T\circ U \circ S^{i-1}\r \cdots \r  T^i\circ U  
\]
as $(\eta^U)^i$.

Associated to $(U,\eta^U)$ there is a functor $U^{\text{gr}}:\Cscr^{\text{gr}}\r
\Dscr^{\text{gr}}$ given by 
\begin{equation}
\label{ref-2-1} 
U^{\text{gr}}(f_i)=(\eta^U)^i_B\circ U(f_i)
\end{equation}
for $f_i\in \Hom^i_\Cscr(A,B)$. It is clear that the formation of $(-)^{\text{gr}}$ is compatible with compositions.

\subsection{Triangulated categories}
We will assume that triangulated categories have a strictly invertible
shift functor. Up to equivalence we may always reduce to this case.
\begin{definitions} An exact  functor $U:\Sscr\r \Tscr$ between triangulated
categories is a graded functor $(U,\eta^U):(\Sscr,[1])\r (\Tscr,[1])$
such that for any distinguished triangle
\[
A\xrightarrow{f} B\xrightarrow{g} C\xrightarrow{h} A[1]
\]
the following triangle
\[
UA\xrightarrow{Uf} UB\xrightarrow{Ug} UC\xrightarrow{\eta^U_A\circ Uh} (UA)[1]
\]
is distinguished.
\end{definitions}
\begin{examples} 
\label{ref-A.3.2-2}
Let $s:\Ascr\r \Ascr$ be the functor which coincides
  with $[1]$ on objects and maps but for which $\eta^s_A:s(A[1])\r (sA)[1]$
  is given by $-\Id_{A[2]}$. Then $(s,\eta^s)$ is an exact endofunctor on $\Ascr$. 
Note in contrast that $[1]$ itself, while being a graded endofunctor, is \emph{not} exact.
\end{examples}
\subsection{Serre functors}
Let $k$ be  a field and assume that $\Cscr$ is a $\Hom$-finite $k$-linear category. 
\begin{definitions} $\Cscr$ satisfies \emph{Serre duality} if there is an
  auto-equivalence $F:\Cscr\r \Cscr$ together with isomorphisms
\begin{equation}
\label{ref-3-3}
\Hom_\Cscr(A,B)\r \Hom_\Cscr(B,FA)^\ast
\end{equation}
natural in $A,B$. Such an $F$ is called a \emph{Serre functor} for $\Cscr$.
\end{definitions}
Putting $B=A$ in \eqref{ref-3-3} yields a canonical element
\[
\Tr_A:\Hom_\Cscr(A,FA)\r k
\]
corresponding to the identity in $\Hom_\Cscr(A,A)$.  It
is easy to see that $\Tr_A(-\circ -)$ defines a non-degenerate pairing
\[
 \Hom_\Cscr(B,FA)\times \Hom_\Cscr(A,B)\r k
\]
and
that the map \eqref{ref-3-3} is given by $f\mapsto: \Tr_A(-\circ f)$.
In addition we have the following fundamental identity \cite{ReVdB1}
\begin{equation}
\label{ref-4-4}
\Tr_A(g\circ f)=\Tr_B(Ff\circ g)
\end{equation}

Now assume that $(\Cscr,S)$ is graded and assume that $\Cscr$ has a Serre
functor $F$. We may make $F$ into a graded functor as follows: we
have to give maps
\[
\eta^F_A:(F\circ S)(A)\r (S\circ F)(A)
\]
natural in $A$.
Using non-degeneracy of the trace pairing we define these
maps via the requirement
\begin{equation}
\label{ref-5-5}
\Tr_A(S^{-1}(\eta^F_A \circ f))=-\Tr_{SA}(f )
\end{equation}
for any $f:SA\r (F\circ S)(A)$. 
\begin{remarks}
The minus sign in this formula is
an arbitrary choice in the graded context, but it is forced in the 
triangulated context. See the proof of Theorem \ref{ref-A.4.4-6} below.
\end{remarks}
\begin{propositions} (Graded Serre duality) For $f_i\in \Hom^i_\Cscr(A,B)$,
$g_{-i}\in \Hom^{-i}_\Cscr(B,FA)$ we have
\[
\Tr_A(g_{-i}\ast f_i)=(-1)^i\Tr_B(F^{\mathrm{gr}}f_i\ast g_{-i})
\]
\end{propositions}
\begin{proof}
We have
\begin{align*}
\Tr_B(F^{\text{gr}}f_i\ast g_{-i})&=\Tr_B(S^{-i}(F^{\mathrm{gr}}f_i)\circ g_{-i})
& \qquad (\text{by } \eqref{ref-1-0})\\
&=\Tr_B(S^{-i}((\eta^F)^i_B \circ F(f_i)\circ S^ig_{-i}))&\qquad (\text{by } \eqref{ref-2-1})\\
&=(-1)^i \Tr_{S^iB}(F(f_i)\circ S^ig_{-i})&\qquad (\text{by }\eqref{ref-5-5})\\
&=(-1)^i \Tr_{A}(S^ig_{-i}\circ f_i)&\qquad (\text{by }\eqref{ref-4-4})\\
&=(-1)^i \Tr_A(g_{-i}\ast f_i)&\qquad (\text{by }\eqref{ref-1-0})&\qed
\end{align*}
\def\qed{}\end{proof}
Assume  now that $\Ascr$ is a $\Hom$-finite $k$-linear triangulated  category
with a Serre functor $F$. 
\begin{theorems} \cite{Bondal4}
\label{ref-A.4.4-6}
$F$ is an  exact functor when equipped with  the graded structure obtained from
 \eqref{ref-5-5} (with $S=[1]$). 
\end{theorems}
\begin{proof}
This is proved
by Bondal and Kapranov in \cite{Bondal4}. We give a somewhat more direct 
proof which makes the signs evident. 

\medskip

We start with a distinguished triangle.
\[
A\xrightarrow{f} B\xrightarrow{g} C\xrightarrow{h} A[1]
\]
We have to construct a map $\delta$ such that the following diagram is a morphism
of distinguished triangles
\[
\xymatrix{%
FA \ar[r]^{Ff} &FB  \ar[r]^{Fg} & FC \ar[r]^{\eta^F_A \circ Fh} & (FA)[1]\\
FA \ar[r]_{Ff}\ar@{=}[u] &FB \ar@{=}[u] \ar[r]_{\alpha} & X \ar[u]^\delta\ar[r]_{\beta} & (FA)[1]\ar@{=}[u]
}
\]
where $X$ is the cone of $Ff$.

In equations:
\begin{align}
\eta^F_A \circ Fh \circ \delta&=\beta \label{ref-6-7}\\
\delta\circ \alpha&=Fg \label{ref-7-8}
\end{align}
For any $x:A\r X[-1]$ we deduce from \eqref{ref-6-7}
\[
(\eta^F_A \circ Fh \circ \delta)[-1]\circ x=\beta[-1] \circ x
\]
Using \eqref{ref-5-5}
this is equivalent to 
\[
\Tr_{A[1]}(Fh \circ \delta \circ x[1])=-\Tr_A(\beta[-1] \circ x)
\]
which using \eqref{ref-4-4} can be further rewritten as
\begin{equation}
\label{ref-8-9}
\Tr_{C}(\delta \circ x[1]\circ h)=-\Tr_A(\beta[-1] \circ x)
\end{equation}
Using the non-degeneracy of the trace pairing we see that \eqref{ref-6-7} is equivalent
to the validity of \eqref{ref-8-9} for all $x:A\r X[-1]$.
Similarly \eqref{ref-7-8} is equivalent to 
\[
\Tr_C(\delta\circ \alpha\circ y)=\Tr_C(Fg\circ y)=\Tr_B(y\circ g)
\]
for all $y:C\r FB$.

Summarizing: we have to find $\delta$ such that the following equations
\begin{equation}
\label{ref-9-10}
\begin{split}
\Tr_{C}(\delta \circ x[1]\circ h)&=-\Tr_A(\beta[-1] \circ x)\\
\Tr_C(\delta\circ \alpha\circ y)&=\Tr_B(y\circ g)
\end{split}
\end{equation}
hold for all $x\in \Hom_\Ascr(A,X[-1])$ and $y\in \Hom_\Ascr(C,FB)$.

We may view the equations \eqref{ref-9-10} as fixing the value of the
function $\Tr_{C}(\delta \circ -)$ on two sub vector spaces of $\Hom_\Ascr(C,X)$. Since
$\Tr_C$ is non-degenerate such a system can be solved provided
 we give the same value on the intersection. Thus we have to show
\[
\alpha\circ y=x[1]\circ h\text{ then }\Tr_B(y\circ g)=-\Tr_A(\beta[-1]\circ x)
\]
To prove this assume $\alpha\circ y=x[1]\circ h$ and consider the following commutative diagram
\begin{equation}
\label{ref-10-11}
\xymatrix{%
FA \ar[r]^{Ff} &FB  \ar[r]^{\alpha} & X \ar[r]^{\beta} & (FA)[1]\\
B\ar@{.>}[u]^\psi\ar[r]_g&C\ar[u]^y\ar[r]_h & A[1] \ar[u]^{x[1]}\ar[r]_{-f[1]}& B[1]\ar@{.>}[u]_{\psi[1]}
}
\end{equation}
where $\psi$ exists because of the axioms of triangulated categories.

\medskip

We compute
\begin{align*}
\Tr_B(y\circ g)&=\Tr_B(Ff\circ \psi)\\
&=\Tr_A(\psi\circ f)\\
&=-\Tr_A(\beta[-1]\circ x)
\end{align*}
In the third line we have used the commutativity of the rightmost square in \eqref{ref-10-11}.
\end{proof}
\subsection{The Calabi-Yau case}
\begin{definitions} A triangulated category with Serre functor $F$ is
  \emph{Calabi-Yau of dimension $n$} if $F\cong s^n$ as graded
  functors, where $s$ is as in Example \ref{ref-A.3.2-2}.
\end{definitions}
\begin{propositions}\label{appendix}
\label{ref-A.5.2-12} Assume that $\Ascr$ is Calabi-Yau of dimension $n$. Then
for $f_i\in \Hom^i_\Ascr(A,B)$, $g_{n-i}\in \Hom^{n-i}_\Ascr(B,A)$ we have
\begin{equation}
\label{ref-11-13}
\Tr_A(g_{n-i}\ast f_i)=(-1)^{i(n-i)} \Tr_B(f_i \ast g_{n-i})
\end{equation}
\end{propositions}
\begin{proof}  We view $g_{n-i}$ as an element of $\Hom_\Ascr^{-i}(A,FB)$
by using the naive identification on objects $(FB)[-i]=B[n][-i]=B[n-i]$.
To avoid confusion we put $h_{-i}=g_{n-i}$.

Graded Serre duality now reads as 
\[
\Tr_A(h_{-i}\ast f_i)=(-1)^i \Tr_B((s^{\text{gr}})^n(f_i)\ast h_{-i})
\]
Writing out everything explicitly we get
\begin{align*}
\Tr_A(h_{-i}[i]\circ f_i)&=(-1)^i \Tr_B\bigl(\bigl((\eta^s)^{ni}_A \circ f_i[n]\bigr)[-i]\circ h_{-i}\bigr)\\
&=(-1)^i \Tr_B\bigl((\eta^s)^{ni}_A[-i] \circ f_i[n-i]\circ h_{-i}\bigr)
\end{align*}
Now composing with $(\eta^s)^{ni}_A[-i]$ is just multiplying by $(-)^{ni}$.
Thus we obtain
\[
\Tr_A(h_{-i}[i]\circ f_i)=(-1)^{i+ni}\Tr_B(f_i[n-i]\circ h_{-i})
\]
which translates into \eqref{ref-11-13}.
\end{proof}
\def\cprime{$'$} \def\cprime{$'$}
\ifx\undefined\bysame
\newcommand{\bysame}{\leavevmode\hbox to3em{\hrulefill}\,}
\fi
\bibliographystyle{amsplain}
\bibliography{calab}

\providecommand{\bysame}{\leavevmode\hbox to3em{\hrulefill}\thinspace}
\providecommand{\MR}{\relax\ifhmode\unskip\space\fi MR }
\providecommand{\MRhref}[2]{%
  \href{http://www.ams.org/mathscinet-getitem?mr=#1}{#2}
}
\providecommand{\href}[2]{#2}
\begin{thebibliography}{10}

\bibitem{RG}
\emph{No references available yet}.

\bibitem{AS}
M.~Artin and W.~Schelter, \emph{Graded algebras of global dimension 3}, Adv.
  Math. \textbf{66} (1987), 171--216.

\bibitem{bergman}
George~M. Bergman, \emph{The diamond lemma for ring theory}, Adv. in Math.
  \textbf{29} (1978), no.~2, 178--218.

\bibitem{BLS}
Raf Bocklandt, Lieven Le~Bruyn, and Stijn Symens, \emph{Isolated singularities,
  smooth orders, and {A}uslander regularity}, Comm. Algebra \textbf{31} (2003),
  no.~12, 6019--6036.

\bibitem{Bondal4}
A.~I. Bondal and M.~M. Kapranov, \emph{Representable functors, {S}erre
  functors, and reconstructions}, Izv. Akad. Nauk SSSR Ser. Mat. \textbf{53}
  (1989), no.~6, 1183--1205, 1337.

\bibitem{kingbutler}
M.~C.~R. Butler and A.~D. King, \emph{Minimal resolutions of algebras}, J.
  Algebra \textbf{212} (1999), no.~1, 323--362.

\bibitem{CBHol}
William Crawley-Boevey and Martin~P. Holland, \emph{Noncommutative deformations
  of {K}leinian singularities}, Duke Math. J. \textbf{92} (1998), no.~3,
  605--635. \MR{MR1620538 (99f:14003)}

\bibitem{oberwolf}
M.~Van den Bergh, \emph{Introduction to super potentials}, M.F.O. Report
  \textbf{6} (2005), 394--396,
  \verb+(http://www.mfo.de/programme/schedule/2005/06/OWR_2005_06.pdf)+.

\bibitem{gap}
The GAP~Group, \emph{{GAP -- Groups, Algorithms, and Programming, Version
  4.4}}, 2005, \verb+(http://www.gap-system.org)+.

\bibitem{orlov}
A.~N. Kapustin and D.~O. Orlov, \emph{Lectures on mirror symmetry, derived
  categories, and {D}-branes}, Uspekhi Mat. Nauk \textbf{59} (2004),
  no.~5(359), 101--134.

\bibitem{Keller}
B.~Keller and D.~Vossieck, \emph{Sous les cat\'egories d\'eriv\'ees}, C. R.
  Acad. Sci. Paris \textbf{305} (1987), 225--228.

\bibitem{groeb}
Teo Mora, \emph{An introduction to commutative and noncommutative {G}r\"obner
  bases}, Theoret. Comput. Sci. \textbf{134} (1994), no.~1, 131--173, Second
  International Colloquium on Words, Languages and Combinatorics (Kyoto, 1992).

\bibitem{ReVdB1}
I.~Reiten and M.~Van~den Bergh, \emph{Noetherian hereditary abelian categories
  satisfying {S}erre duality}, J. Amer. Math. Soc. \textbf{15} (2002), no.~2,
  295--366 (electronic).

\bibitem{Verdier}
J.-L. Verdier, \emph{Des cat\'egories dériv\'ees des cat\'egories
  ab\'eliennes}, Ast\'erisque, S.M.F.

\end{thebibliography}

\end{document}